\documentclass[12pt]{amsart}

\usepackage{pb-diagram}

\newtheorem{thm}{Theorem}[section]    
\newtheorem{cor}[thm]{Corollary}
\newtheorem{lem}[thm]{Lemma}
\newtheorem{prop}[thm]{Proposition}

\theoremstyle{remark}                
\newtheorem{rmk}[thm]{Remark}
\newtheorem{ex}[thm]{Example}

\theoremstyle{definition}

\numberwithin{equation}{section}

\newcommand{\C}{\mathbb C}
\newcommand{\csa}{\ensuremath{C^*}-algebra}
\newcommand{\ib}{imprimitivity bimodule}

\newcommand{\secref}[1]{Section~\ref{#1}}
\newcommand{\thmref}[1]{Theorem~\ref{#1}}

\newcommand{\propref}[1]{Proposition~\ref{#1}}
\newcommand{\lemref}[1]{Lemma~\ref{#1}}
\newcommand{\eqnref}[1]{\eqref{#1}}

\newcommand{\im}{\operatorname{im}}
\newcommand{\Aut}{\operatorname{Aut}}
\newcommand{\Ind}{\operatorname{Ind}}

\newcommand{\Res}{\operatorname{Res}}
\newcommand{\Rep}{\operatorname{Rep}}
\newcommand{\clsp}{\overline{\operatorname{span}}}
\newcommand{\clspn}{\clsp}
\newcommand{\supp}{\operatorname{supp}}
\newcommand{\id}{\operatorname{id}}
\newcommand{\spn}{\operatorname{span}}
\newcommand{\Ad}{\operatorname{Ad}}

\newcommand{\x}{\times}
\newcommand{\ox}{\otimes}
\newcommand{\rx}{\times}
\newcommand{\comp}{\circ}

\newcommand{\inv}[1]{{#1}^{-1}}

\newcommand{\slashstyle}{\mathcal}

\renewcommand{\H}{{\slashstyle H}}
\renewcommand{\L}{{\slashstyle L}}

\newcommand{\K}{{\slashstyle K}}

\newcommand{\LtwoG}{{L^2(G)}}

\newcommand{\lip}[3]{
  {\vphantom\langle}_{#1}\!\!\left\langle{#2},{#3}\right\rangle  }
\newcommand{\rip}[3]{
  \left\langle{#2},{#3}\right\rangle_{\!{#1}}  }

\newcommand{\at}{{\Tilde{\alpha}}}
\newcommand{\ah}{{\widehat{\alpha}}}
\newcommand{\ahh}{{\widehat{\widehat{\alpha}}}}   

\newcommand{\AxaH}{A\rx_\alpha H}

\newcommand{\AxaG}{A\rx_\alpha G}
\newcommand{\AxG}{A\rx_\alpha G}

\newcommand{\AxGxG}{(\AxaG)\rx_{\ah} G}

\newcommand{\CGAxG}{C_0(G,A)\rx_{\alpha\ox\tau}G}
\newcommand{\CGmHAxG}{C_0(G/H,A)\rx_{\alpha\ox\tau}G}

\newcommand{\CGAxGxH}{(\CGAxG)\rx_{\ahh} H}
\newcommand{\CGAxHxG}{C_0(G,A)\rx_{\at}(H\x G)}
\newcommand{\CGmNAxG}{C_0(G/N,A)\rx_{\alpha\ox\tau}G}

\newcommand{\D}{{\mathcal D}}                    
\newcommand{\I}{{\mathcal I}}                    

\newcommand{\dx}{\delta_X}
\newcommand{\da}{\delta_A}
\newcommand{\db}{\delta_B}
\newcommand{\dl}{\delta_L}
\newcommand{\dg}{\delta_G}

\begin{document}

\title[Crossed products by dual coactions]{Crossed products by dual coactions
of groups and homogeneous spaces}

\date{20 August, 1996}

\author[Echterhoff]{Siegfried Echterhoff}
\address{
Fachbereich Mathematik-Informatik\\
Universit\"at-Gesamthochschule Paderborn\\
D-33095 Paderborn\\ Germany}
\email{echter@uni-paderborn.de}

\author[Kaliszewski]{S.~Kaliszewski}
\address{Department of Mathematics \\ University of Newcastle \\
NSW 2308 \\ Australia}
\email{kaz@frey.newcastle.edu.au}

\author[Raeburn]{Iain Raeburn}
\address{Department of Mathematics \\ University of Newcastle \\
NSW 2308 \\ Australia}
\email{iain@frey.newcastle.edu.au}

\subjclass{46L55, secondary 22D25}

\keywords{$C^*$-algebra, coaction, crossed product,
imprimitivity, homogeneous space}
\thanks{This research was supported by the Australian Research Council.}

\begin{abstract}
Mansfield showed how to induce representations of crossed products of
$C^*$-algebras by coactions from crossed products by quotient groups and proved
an imprimitivity theorem characterising these induced representations. 
We give an alternative construction of his
bimodule in the case of dual coactions, based on the symmetric
imprimitivity theorem of the third author; this
provides a more workable way of inducing representations of crossed products of
$C^*$-algebras by dual coactions.  The construction works for
homogeneous spaces as well as quotient groups, and we prove an imprimitivity
theorem for these induced representations.
\end{abstract}

\maketitle


Coactions of groups on $C^*$-algebras, and their crossed
products, were introduced to make duality arguments available
for the study of dynamical systems involving actions of
nonabelian groups. For these to be effective, one needs to
understand the representation theory of crossed products by
coactions. The most powerful tool we have was provided by Mansfield
\cite{ManJF}: he showed how to induce representations from
crossed products by quotient groups, and proved an
imprimitivity theorem which characterises these induced
representations. Unfortunately, Mansfield's construction is
complicated and technical. The Hilbert bimodule with which he
defines induced representations is difficult to manipulate,
and one is tempted to seek other realisations of this bimodule
and the induced representations. Here we show that, at
least for the dual coactions arising in the study of ordinary 
dynamical systems, there is an alternative bimodule  built
along more conventional lines from spaces of continuous
functions with values in $C^*$-algebras. This bimodule will 
be easier to work with, and will allow us to induce
representations from quotient homogeneous spaces as well as
quotient groups. 

The core of our construction is a special case of the
symmetric imprimitivity theorem of \cite{RaeIC}. Suppose
$\alpha$ is an action of a locally compact group $G$ on a
$C^*$-algebra $A$. For each closed subgroup $H$ of $G$, there
is a diagonal action $\alpha\otimes\tau$ of $G$ on
$A\otimes C_0(G/H)$: if we identify $A\otimes
C_0(G/H)$ with $C_0(G/H,A)$ in the usual way, then
$(\alpha\otimes\tau)_t(f)(sH)=f(t^{-1}sH)$. We show in
\S1 that there is a natural Morita equivalence between an
iterated crossed product
$(C_0(G,A)\times_{\alpha\otimes\tau}G)
\times H$ and the
imprimitivity algebra $C_0(G/H,A)\times_{\alpha\otimes\tau}G$
of Green \cite{GreLS}. If $H$ is normal, this imprimitivity
algebra can be identified with the crossed product
$(A\times_\alpha G)\times_{\widehat\alpha|}G/H$ by the
restriction of the dual coaction, and the iterated crossed
product with $((A\times_\alpha G)\times_{\widehat\alpha}
G)\times_{\widehat{\widehat\alpha}|}H$; the existence of our Morita
equivalence is therefore predicted by Mansfield's
imprimitivity theorem, although his construction gives no hint
that the bimodule can be realised as a completion of
$C_c(G\times G,A)$. In \S2, we shall discuss these
isomorphisms in detail, and show how our bimodule can be used
to induce representations from $G/H$ to $G$ even when $H$ is
not normal.

Although it is not clear in general how to define 
coactions of homogeneous spaces, let alone their crossed products 
(see the discussion at the start of \S2), there is considerable
evidence that our inducing process is a step in the right direction.
There is an appropriate imprimitivity theorem
(Proposition \ref{imp-homo-sp}), the induction process
interacts  with Green induction and duality as one would expect
from the results of \cite{EchDI} and \cite{KQR-DR} (Theorem
\ref{HGK-thm} and Corollary \ref{ind-res-cor}), and our bimodule
is isomorphic to Mansfield's when the subgroup
$H$ is normal and  amenable (Theorem \ref{DK-thm}).

When the subgroup $H$ is normal but not amenable, the
relationship between our bimodule and the extension of
Manfield's in \cite[\S3]{KQ-IC} becomes quite subtle. There are
two candidates for the crossed product $(A\times G)\times
G/H$: the spatial version on $\H\otimes L^2(G)$ used in
\cite{KQ-IC}, and the imprimitivity algebra $C_0(G/H,A)\times G$.
We believe that one can usefully view the former as a reduced
crossed product by the homogeneous space, and the latter as a
full crossed product. We discuss this in detail in
\S2. However, that the two can be different has an interesting
consequence: the bimodule used in \cite{KQ-IC} can be a proper quotient of
the one we construct in \S1. Thus for nonamenable subgroups, our Morita
equivalence is analogous to Green's equivalence of $A\times H$
and $C_0(G/H,A)\times G$, whereas  \cite[Theorem
3.3]{KQ-IC} is analogous to that of the reduced crossed products 
$A\times_r H$
and $C_0(G/H,A)\times_r G$.

While we are discussing crossed products by homogeneous
spaces, it is worth pointing out that for any
coaction $(B,G,\delta)$ and any closed subgroup $H$, the
spatially defined algebra
$B\times G/H$ is Morita equivalent to 
$(B\times_\delta G)\times_{\widehat\delta,r}H$; however, this
equivalence is obtained as a composition of other
equivalences, and is not obviously implemented by any one
concretely defined bimodule. We discuss this weak version of
Mansfield's Imprimitivity Theorem in an appendix.
 
\section*{ Preliminaries}

Let $G$ be a locally compact group; we always use left Haar measure on
$G$. We  denote by
$\lambda$ the left regular representation of $G$ on $L^2(G)$, and by $M$
the representation of $C_0(G)$ by multiplication operators on $L^2(G)$.
We extend representations and nondegenerate homomorphisms to multiplier
algebras without comment or change of notation; thus, for example, $M$
also denotes the representation of $C_b(G)=M(C_0(G))$ by multiplication
operators.

An action of $G$ on a $C^*$-algebra $A$ is a
homomorphism $\alpha$ of  $G$ into $\Aut A$ such that $s\mapsto
\alpha_s(a)$ is continuous for every $a\in A$.   The crossed product
$(\AxG,i_A,i_G)$ is the universal object for covariant representations
of $(A,G,\alpha)$, as in \cite{RaeOC}; the set $C_c(G,A)$ of continuous,
compactly supported functions from $G$ into $A$ embeds as a dense
*-subalgebra of $\AxG$, with 
\[
f*g(s) = \int_G f(t)\alpha_t(g(t^{-1}s))\, dt\quad\text{and}\quad
f^*(s) = \alpha_s(f(s^{-1})^*)\Delta_G(s)^{-1}.
\]
If $\pi$ is a nondegenerate representation of $A$, the induced
representation $\Ind\pi$ of the system $(A,G,\alpha)$ is the covariant
representation $(\widetilde\pi,1\otimes\lambda)$, in which
$\widetilde\pi(a)\xi(s):=\pi(\alpha_s^{-1}(a))(\xi(s))$ for $\xi\in
L^2(G,\H_\pi)=\H_\pi\otimes L^2(G)$.  If $H$ is a closed subgroup of $G$,
we identify $A\otimes C_0(G/H)$ with $C_0(G/H,A)$; we write
$\alpha\otimes\tau$ for the diagonal action of $G$ on either, so that
$(\alpha\otimes\tau)_t(f)(sH)=f(t^{-1}sH)$ for $f\in C_0(G/H,A)$.
We use $\sigma$ to denote the action of $G$ on $C_0(G)$ by right
translation: $\sigma_t(f)(s) := f(st)$.  

We use the full coactions of \cite{RaeCR}, as modified in \cite{QuiFR}: we
use minimal tensor products throughout. Thus a coaction
$\delta$ of $G$ on a $C^*$-algebra $B$ is a nondegenerate homomorphism
$\delta\colon B\to M(B\otimes C^*(G))$ such that 
\begin{equation*}
(\delta\otimes\id)\circ\delta = (\id\otimes\delta_G)\circ\delta
\quad\text{and}\quad
\delta(b)(1\otimes z)\in B\otimes C^*(G)
\end{equation*}
for all $b\in B$ and $z\in C^*(G)$,
where $\delta_G\colon C^*(G)\to M(C^*(G)\otimes C^*(G))$ is the
comultiplication on $C^*(G)$ characterised by $\delta_G(i_G(s))=
i_G(s)\otimes i_G(s)$.   If $N$ is a closed
normal subgroup of $G$ and $q\colon C^*(G)\to M(C^*(G/N))$ is
characterised by $q(i_G(s))=i_{G/N}(sN)$, then
$(\id\otimes q)\circ\delta$ is a coaction of $G/N$ on $B$, called the
restriction of $\delta$ to $G/N$, and denoted $\delta|$.  The crossed
product
$(B\times_\delta G,j_B,j_{C(G)})$ is the universal object for covariant
representations of
$(B,G,\delta)$; in particular
\[
B\times_\delta G = \clspn\{ j_B(b)j_{C(G)}(f)\mid b\in B, f\in
C_0(G)\}.
\]
If $\pi$ is a nondegenerate representation of $B$, the induced
representation $\Ind\pi$ of $(B,G,\delta)$ is the covariant
representation $((\pi\otimes\lambda)\circ\delta,1\otimes M)$ on
$\H_\pi\otimes L^2(G)$. We shall follow the conventions of \cite{RaeCR}
concerning dual actions and coactions.


\section{The symmetric imprimitivity theorem}\label{closed-sec}

We begin by recalling the symmetric
imprimitivity theorem of \cite{RaeIC}.  Our conventions will be slightly
different from those used there; here 
the second group $L$ acts on the
right of the locally compact space $P$. To convert to the
two-left-actions situation of \cite{RaeIC}, just let 
$l\cdot p = p \cdot\inv{l}$.  

Consider a \csa\ $D$, two locally compact groups $K$ and $L$, and a
locally compact space $P$; suppose that $K$ acts freely and properly on
the left of $P$, and that $L$ acts likewise on the right, and that 
these actions commute (i.e. $k\cdot(p\cdot l) = (k\cdot p)\cdot l)$.  
Suppose also that we have commuting actions $\sigma$  of 
$K$ and $\rho$ of $L$ on $D$. 
Recall that for the left action of $K$ we define
the induced \csa\ $\Ind\sigma$ to be the set of continuous bounded
functions $f\colon P\to D$ such that
$f(k\cdot p) = \sigma_k(f(p))$
for all $k\in K$ and $p\in P$, and such that the function 
$Kp\mapsto\|f(p)\|$ vanishes at infinity on $K\setminus P$.  For the
right action of $L$ we define the induced \csa\ $\Ind\rho$ to be the
set of continuous bounded functions $f\colon P\to D$ such that
$f(p\cdot l) = \inv{\rho_l}(f(p))$
for all $p\in P$ and $l\in L$, and such that the function 
$pL\mapsto \|f(p)\|$ vanishes at infinity on $P/L$.  The induced
algebras are \csa s with pointwise operations, and
carry actions $\gamma\colon K\to \Aut(\Ind\rho)$ and
$\delta\colon L\to \Aut(\Ind\sigma)$ given by 
$$
\gamma_k(f)(p) = \sigma_k(f(\inv{k}\cdot p))\quad\text{and}\quad
\delta_l(f)(p) = \rho_l(f(p\cdot l)).
$$
Then \cite[Theorem~1.1]{RaeIC}\ states that $C_c(P,D)$
can be given a pre-\ib\ structure which
completes to give a Morita equivalence between
$\Ind\rho\rx_\gamma K$ and $\Ind\sigma\rx_\delta L$.
The actions and inner products are given for 
$b\in C_c(K,\Ind\rho)\subseteq\Ind\rho\rx_\gamma K$, $x$ and $y$ in 
$C_c(P,D)$, and $c\in C_c(L,\Ind\sigma)\subseteq\Ind\sigma\rx_\delta L$
as follows:
\begin{eqnarray*}
b\cdot x(p) &=& \int_K b(t,p) \sigma_t(x(\inv{t}\cdot p))\,
\Delta_K(t)^{\frac12}\, dt\\
x\cdot c(p) &=& \int_L \rho_s\bigl(x(p\cdot s) c(\inv{s},p\cdot s)\bigr)\,
\Delta_L(s)^{-{\frac12}}\, ds\\
\lip{\Ind\rho\rx_\gamma K}{x}{y}(k,p) &=& 
\int_L \rho_s\bigl(x(p\cdot s)
\sigma_k(y(\inv{k}\cdot p\cdot s)^*\bigr)\, ds\,
\Delta_K(k)^{-{\frac12}}\\
\rip{\Ind\sigma\rx_\delta L}{x}{y}(l,p) &=&
\int_K \sigma_t\bigl(x(\inv{t}\cdot p)^*
\rho_l(y(\inv{t}\cdot p\cdot l))\bigr)\,
dt\, \Delta_L(l)^{-{\frac12}}.
\end{eqnarray*}

If $\alpha:G\to \Aut A$ is an action, we denote by 
$\widehat{\widehat{\alpha}}$ the action of $G$ on 
$C_0(G,A)\times_{\alpha\otimes\tau}G$
given for $f\in C_c(G\times G,A)$ by
$$
\ahh_t(f)(r,s) = f(r,st).
$$
(This action is carried into the usual second dual action on 
$\AxGxG\cong C_0(G,A)\times_{\alpha\otimes \tau}G$ under the isomorphism of
Lemma \ref{fund-lem} below.)

\begin{prop}\label{R-bimod-prop}
Let $\alpha:G\to \Aut A$ be an action, and let $H$ be a closed 
subgroup of $G$.  Then there exists a pre-imprimitivity bimodule
structure on $C_c(G\times G,A)$ which completes to give an 
$(C_0(G,A)\times_{\alpha\otimes \tau} G)
\times_{\ahh|}H$ -- $C_0(G/H,A)\times_{\alpha\otimes \tau}G$ imprimitivity
bimodule.
\end{prop}

\begin{proof}
We apply the symmetric imprimitivity theorem, with 
$P=G\x G$, $K = H\x G$, $L=G$, and $D=A$.  Define a left action of 
$H\x G$, and a right action of $G$ on $G\x G$ by
\begin{equation}\label{P-act-eqns}
(h,t)\cdot(r,s) = (hr,ts)\quad\text{and}\quad
(r,s)\cdot t = (rt,st).
\end{equation}
Both these actions are free and proper, and they
commute with one another.  Define actions $\sigma$ and $\rho$ of 
$H\x G$ and $G$ on $A$ as follows:
$$
\sigma_{(h,t)}(a) = \alpha_t(a)\quad\text{and}\quad
\rho_t(a) = a.
$$
It is clear that these actions also commute; thus by the symmetric
imprimitivity theorem \cite[Theorem~1.1]{RaeIC}, 
$C_c(G\x G,A)$ completes to give a
$
\Ind\rho\rx_\gamma(H\x G)$ --  $\Ind\sigma\rx_\delta G.
$ imprimitivity bimodule.

It only remains to identify $\Ind\rho\rx_\gamma(H\x G)$ with 
$(C_0(G,A)\times_{\alpha\otimes\tau}G)\times_{\ahh}H$,
and $\Ind\sigma\rx_\delta G$ with
$C_0(G/H,A)\times_{\alpha\otimes\tau}G$.
To this end, we first remark that $H\x G$ acts on $C_0(G,A)$ by
$
\at_{(h,t)}(f)(s) = \alpha_t(f(\inv{t}sh)),
$
and then that the identity map of $C_c(H\x G\x G,A)$ onto itself
extends to an isomorphism of $\CGAxHxG$ onto $\CGAxGxH$.

Next, note that $(G\x G)/G$ (with the action of \eqnref{P-act-eqns}) is
homeomorphic to $G$ via the map $(r,s)\mapsto s\inv{r}$, so we have a
bijection $\Theta\colon C_0(G,A)\to\Ind\rho$ given by
$$
\Theta(f)(r,s) = f(s\inv{r}),\quad
\inv{\Theta}(g)(s) = g(e,s).
$$
Since the operations on both $C_0(G,A)$ and $\Ind\rho$ are pointwise, 
$\Theta$ gives an isomorphism of the \csa s.  

Now $\Theta$ is $\at$ -- $\gamma$ equivariant:
\begin{eqnarray*}
\Theta(\at_{(h,t)}(f))(r,s) & = & \at_{(h,t)}(f)(s\inv{r})\\
 & = & \alpha_t(f(\inv{t}s\inv{r}h))\\
 & = & \alpha_t(\Theta(f)(\inv{h}r,\inv{t}s))\\
 & = & \sigma_{(h,t)}\left(\Theta(f(\inv{(h,t)}\cdot(r,s)))\right)\\
 & = & \gamma_{(h,t)}(\Theta(f))(r,s).
\end{eqnarray*}
Thus $\Theta$ induces an isomorphism of 
$\CGAxHxG$ onto $\Ind\rho\times_\gamma(H\times G)$.  
Combined with the previous isomorphism,
this completes the first identification.  

For the second identification, note that $(H\x G)\setminus (G\x G)$ (with
the action of \eqnref{P-act-eqns}) is homeomorphic to $G/H$ via the map
$(r,s)\mapsto \inv{r}H$.  So we have a bijection $\Omega\colon
C_0(G/H,A)\to \Ind\sigma$ given by 
$$
\Omega(f)(r,s) = \alpha_s(f(\inv{r}H)),\quad
\inv{\Omega}(g)(tH) = g(\inv{t},e).
$$

As above, since the operations on both algebras are pointwise,
$\Omega$ is an isomorphism.  Moreover, $\Omega$ is
$\alpha\ox\tau$ -- $\delta$ equivariant:
\begin{eqnarray*}
\Omega(\alpha_t\ox\tau_t(f))(r,s) & = &
\alpha_s(\alpha_t\ox\tau_t(f)(\inv{r}H))\\
 & = & \alpha_s(\alpha_t(f(\inv{t}\inv{r}H)))\\
 & = & \alpha_{st}(f(\inv{(rt)}H))\\
 & = & \Omega(f)(rt,st)\\
 & = & \delta_t(\Omega(f))(r,s).
\end{eqnarray*}
Thus $\Omega$ 
induces the second identification of crossed products.
\end{proof}

The isomorphisms of the proof of
\propref{R-bimod-prop}\  can be used to
make $C_c(G\x G,A)$ explicitly a
$C_c(H\times G\times G,A)$ -- $C_c(G\times G/H, A)$
pre-\ib. However, for technical reasons, 
we shall  combine these with the automorphism
$\Upsilon$ of $C_c(G\times G,A)$ defined by
\[
\Upsilon(x)(r,s) = x(r,rs^{-1})\Delta_G(r)^{\frac12}.
\]
This gives a  bimodule structure which is more natural for our
considerations in
\secref{mansfield-sec}.  The resulting  
actions and inner products are given, for $f\in
C_c(H\x G\x G,A)$, $x$ and $y$ in $C_c(G\times G,A)$,
and $g\in C_c(G\x G/H,A)$ as follows:
\begin{eqnarray}
f\cdot x(r,s) &=& \int_G\int_H
f(h,t,s)\alpha_t(x(\inv{t}r,\inv{t}sh))\,
\Delta_H(h)^{\frac12}\,
dh\, dt\label{K-eqn-a}\\
x\cdot g(r,s) &=& \int_G x(t,s)
\alpha_{t}(g(\inv{t}r,\inv{t}sH))\, dt
\label{K-eqn-b}\\
\lip{L}{x}{y}(h,r,s) &=& \int_G x(t,s)
\alpha_r(y(\inv{r}t,\inv{r}sh)^*)\, \Delta_H(h)^{-{\frac12}}\,
\Delta_G(r^{-1}t)\, dt\,
\label{K-eqn-c}\\
\rip{R}{x}{y}(r,sH) &=& \int_G \int_H
\alpha_t(x(\inv{t},\inv{t}sh)^* y(\inv{t}r,\inv{t}sh))\,
\Delta_G(t^{-1})\, dh\, dt\label{K-eqn-d}.
\end{eqnarray}


\section{Inducing representations from homogeneous spaces}

It is a major defect of the current theory of crossed
products by coactions that we do not know how to best define
coactions of homogeneous spaces and their crossed products.
However, if we start with a coaction of $G$ on $B$, and $H$
is a closed subgroup of $G$, we can obtain what should be
covariant  representations of $(B,G/H,\delta)$ by restricting
covariant representations $(\pi,\mu)$ of $(B,G,\delta)$: 
just extend $\mu$ to the multiplier algebra
$M(C_0(G))=C_b(G)$ and restrict it to the subalgebra
$C_0(G/H)$ of functions constant on $H$-cosets. In
particular, we can restrict a regular representation
$((\pi\otimes\lambda)\circ\delta,1\otimes M)$, and define the
\emph{reduced crossed product} $B\times_{\delta,r}G/H$ to be the
$C^*$-subalgebra of
$B(\H_\pi\otimes L^2(G))$ generated by the operators
\[
\{(\pi\otimes\lambda)\circ\delta(b)(1\otimes M_f) \mid b\in
B,f\in C_0(G/H)\}.
\]
Provided $\Ind\pi$ is faithful on $B\times_\delta G$,
$M(B\times_\delta G)$ is represented faithfully on
$\H_\pi\otimes L^2(G)$, so $B\times_{\delta,r}G/H$ is actually
a subalgebra of $M(B\times_\delta G)$; thus, with the proviso that
$\Ind\pi$ is faithful, the isomorphism class of
$B\times_{\delta,r}G/H$ is independent of the choice of
$\pi$.

\begin{rmk}  (1) We have chosen the notation
$B\times_{\delta,r}G/H$ to stress that the reduced crossed
product depends on the coaction $\delta$, and, implicitly, on
the group $G$. (A given space may be realisable in several
different ways as a homogeneous space.) For a normal subgroup
$N$, $B\times_{\delta,r}G/N$ is not necessarily the same as the crossed
product
$B\times_{\delta|}G/N$.
Restricting the regular representation of $(B,G,\delta)$
gives a covariant representation of $(B,G/N,\delta|)$ on
$\H_\pi\otimes L^2(G)$, which is known to be faithful if $N$
is amenable (\cite[Lemma 3.2]{KQ-IC}, or Corollary \ref{reg-rep-faith}
below), but is not faithful in general (Remark \ref{not-faith} below).

The notation we have chosen is consistent with that used by
Mansfield to distinguish the subalgebra of $B(\H_\pi\otimes
L^2(G))$ from his spatially defined crossed product
$B\times_{\delta|}G/N$ on
$\H_\pi\otimes L^2(G/N)$. We mention in passing that, for
arbitrary $H$, it follows from \cite[Proposition 8]{ManJF} that
\[
B\times_{\delta,r}G/H=
\clsp\{(\pi\otimes\lambda)\circ\delta(b)(1\otimes M_f)\mid b\in
B,f\in C_0(G/H)\}.
\]

(2) Since $M(B\times_\delta G)$ is faithfully represented on
$\H_\pi\otimes L^2(G)$, for normal $N$ the algebra
$B\times_{\delta,r}G/N$ is the algebra $\im(j_B\times j_G|)$
appearing in \cite[Theorem 3.3]{KQ-IC}, and hence that theorem establishes
a Morita equivalence between $B\times_{\delta,r}G/N$ and the
reduced crossed product $(B\times_\delta G)\times_r N$. This
bimodule can be used to define induction of representations
from $B\times_{\delta,r}G/N$ to $B\times_\delta G$. As we
shall see,  this is not necessarily 
the same as the induction process we shall construct for $B$
of the form $A\times_\alpha G$.

\end{rmk}

When $\delta$ is the dual coaction $\widehat\alpha$ of an action
$\alpha:G\to\Aut A$, there is also a natural candidate for a
full crossed product
$B\times_\delta G/H$, whose representations are given by
certain covariant pairs
$(\pi,\mu)$ of representations of $B$ and $C_0(G/H)$. To
motivate this, we recall that for normal $N$, the crossed
product $(A\times_\alpha G)\times_{\widehat\alpha|}G/N$ is one
realisation of Green's imprimitivity algebra $(A\otimes
C_0(G/N))\times_{\alpha\otimes\tau}G$; indeed, the
resulting interpretation of Green's Imprimitivity Theorem
motivated Mansfield's theorem (see \cite{ManIR}). We digress to
establish this realisation in the context of
full coactions and nonamenable subgroups.

\begin{lem}\label{rae-co-lem}
Let $\alpha:G\to\Aut A$ be an action, and let
$N$ be a closed normal subgroup of
$G$. Consider representations 
$\pi$, $U$, and $\mu$ of $A$,
$G$ and $C_0(G/N)$, respectively, 
on a Hilbert space $\H$.  Then $(\pi,U)$
is a covariant representation of $(A,G,\alpha)$ and $(\pi\x U,
\mu)$ is a covariant representation of $(\AxaG,G/N,\ah)$ if
and only if $\pi$ and $\mu$ have commuting ranges and 
$(\pi\ox\mu,U)$ is a covariant representation of $(C_0(G/N,A),G,
\alpha\ox\tau)$.  
\end{lem}

\begin{proof}
The proof is sketched in \cite[Example~2.9]{RaeCR}. 
\end{proof}

\begin{lem}\label{fund-lem}
Let $\alpha:G\to \Aut A$ be an action of a locally compact group, and
let $N$ be a closed normal subgroup of
$G$. Then there is an isomorphism 
$$\Psi\colon\CGmNAxG\to (A\times_{\alpha}G)\times_{\widehat\alpha|}G/N$$
which is natural in the sense that 
\begin{equation*}
\Psi\comp k_A = j_{A\rx G}\comp
i_A, \quad\Psi\comp k_G = j_{A\rx G}\comp i_G,
\quad\text{and}\quad 
\Psi\comp k_{C(G/N)} = j_{C(G/N)},
\end{equation*}
 where
$(k_A\ox k_{C(G/N)},k_G)$ are the canonical maps of 
$(C_0(G/N,A),G,\alpha\ox\tau)$ into
 the crossed product.  The induced map on
representations takes $(\pi\x U)\x\mu$ to $(\pi\ox\mu)\x U$.  
\end{lem}

\begin{proof}
Realise $\CGmNAxG$ on  $\H$; then
$k_A$ and $k_{C(G/N)}$ are commuting representations on $\H$, and
$(k_A\ox k_{C(G/N)},k_G)$ is a covariant representation of 
$(C_0(G/N,A),G,\alpha\ox\tau)$.  By \lemref{rae-co-lem}, 
$(k_A,k_G)$ is covariant for 
$(A,G,\alpha)$, and
$(k_A\x k_G, k_{C(G/N)})$ is covariant for
$(\AxaG,G/N,\ah|)$.  It follows  that there is a 
nondegenerate representation
$\Phi = (k_A\x k_G)\x k_{C(G/N)}$ of 
$(A\times_{\alpha}G)\times_{\widehat\alpha|}G/N$ on $\H$ such that
\begin{equation}\label{phi-eqn-a}
\Phi\comp j_{A\rx G}\comp i_A = k_A,\quad 
\Phi\comp j_{A\rx G}\comp i_G = k_G,\quad
\Phi\comp j_{C(G/N)} = k_{C(G/N)}.
\end{equation}

Now suppose  
$(A\times_{\alpha}G)\times_{\ah|}G/N$ acts on $\K$.  Then
\begin{equation*}
(j_{\AxG},j_{C(G/N)}) = 
((j_{A\rx G}\comp i_A)\x(j_{A\rx G}\comp i_G), j_{C(G/N)})
\end{equation*}
is a covariant representation of
$(\AxaG,G/N,\ah|)$.  Thus we deduce from \lemref{rae-co-lem} 
that
$((j_{A\rx G}\comp i_A)\ox j_{C(G/N)}, j_{A\rx G}\comp i_G)$
is covariant for $(C_0(G/N,A),G,\alpha\ox\tau)$, and hence there is a
representation 
$\Psi = ((j_{A\rx G}\comp i_A)\ox j_{C(G/N)})\x (j_{A\rx G}\comp i_G)$
of $\CGmNAxG$ on $\K$ such that
\begin{equation}\label{psi-eqn-a}
\Psi\comp k_A = j_{A\rx G}\comp i_A,\quad 
\Psi\comp k_G = j_{A\rx G}\comp i_G,\quad
\Psi\comp k_{C(G/N)} = j_{C(G/N)}.  
\end{equation}
Equations~\eqref{phi-eqn-a}
and  \eqref{psi-eqn-a} imply 
that $\Psi$ is an inverse for $\Phi$.

For the last statement,  let 
$(\pi\x U)\x \mu$ be a 
representation of $(A\times_{\alpha}G)\times_{\ah|}G/N$.  With $a\in A$,
$z\in C_c(G)$, and $f\in C_c(G/N)$, 
$k_A\ox k_{C(G/N)}(a\ox f)k_G(z)$ 
is a typical enough element of $\CGmNAxG$, and  we have:
\begin{eqnarray*}
\lefteqn{\left((\pi\x U)\x\mu\right)
         \comp\Psi\left(k_A\ox k_{C(G/N)}(a\ox
f)k_G(z)\right)}{\mbox{\qquad\qquad\qquad\qquad}}\\
 & = &  \left((\pi\x U)\x\mu\right)
        \left(j_{C(G/N)}(f)j_{A\rx G}(i_A(a)i_G(z))\right)\\
 & = &  \mu(f)\pi(a) U(z)\\
 & = &  \pi\ox\mu(a\ox f) U(z)\\
 & = &  \left((\pi\ox\mu)\x U\right)(k_A\ox k_{C(G/N)}(a\ox f)k_G(z)).
\end{eqnarray*}
\end{proof}

For the rest of this section, $\alpha:G\to\Aut A$ will be an action of
a locally compact group, and $H$ an arbitrary closed subgroup of $G$.
Lemma
\ref{rae-co-lem} suggests that it is reasonable to say that a pair of
representations
$(\pi\times U,\mu)$ of
$(A\times_\alpha G,C_0(G/H))$ is a \emph{covariant representation of
$(A\times_\alpha G,G/H,\widehat\alpha)$} if the ranges of $\pi$
and $\mu$ commute and $(\pi\otimes\mu,U)$ is a covariant
representation of $(C_0(G/H,A),G,\alpha\otimes\tau)$. Then we
can view $C_0(G/H,A)\times_{\alpha\otimes\tau}G$ as a full
crossed product of $A\times_\alpha G$ by the ``coaction''
$\widehat\alpha$ of the homogeneous space $G/H$. We now want to discuss the
``regular representations'' of this full crossed product. But first we
need to know that certain representations of $A\times_\alpha G$ induce
to faithful representations of $(A\times_\alpha
G)\times_{\widehat\alpha} G$, so that we can use them to define the
reduced crossed product $(A\times_\alpha G)\times_{\widehat\alpha,r}G/H$.

\begin{lem}\label{ind-faith}
Let $(\pi,U)$ be a covariant representation of $(A,G,\alpha)$ such that
$\pi$ is faithful. Then the representation $\Ind(\pi\times U)$ of
$(A\times_\alpha G)\times_{\widehat\alpha} G$ is faithful; so is the
corresponding representation $(\pi\otimes M)\times(U\otimes \lambda)$ of
$(A\otimes C_0(G))\times_{\alpha\otimes \tau}G$.
\end{lem}

\begin{proof}
Since 
\begin{eqnarray*}
\Ind(\pi\times U)&=&(((\pi\times
U)\otimes\lambda)\circ\widehat\alpha)\times (1\otimes M)\\
&=&
((\pi\otimes 1)\times(U\otimes\lambda))\times (1\otimes M),
\end{eqnarray*}
it follows from Lemma
\ref{fund-lem} that it is enough to show that the representation
$(\pi\otimes M)\times (U\otimes\lambda)$ of
$C_0(G,A)\times_{\alpha\otimes\tau}G$ is faithful. The automorphism
$\phi$ of $C_0(G,A)$ defined by $\phi(f)(t)=\alpha_{t^{-1}}(f(t))$
induces an isomorphism of $C_0(G,A)\times_{\alpha\otimes\tau}G$ onto
\[
C_0(G,A)\times_{\id\otimes\tau}G\cong A\otimes (C_0(G)\times_\tau G)\cong
A\otimes\K(L^2(G)).
\]
If we now define $V$ on $L^2(G,\H)$ by $V\xi(t)=U_t(\xi(t))$, then one
can verify that
\[
V^*(\pi\otimes M)(\phi^{-1}(f))V=\pi\otimes M(f),\ \text{ and}\ 
V^*(U\otimes\lambda)V=1\otimes\lambda.
\] 
Since the representation $(\pi\otimes M)\times
(1\otimes\lambda)=\pi\otimes (M\times\lambda)$ is certainly faithful
on $A\otimes\K(L^2(G))$,
the result follows.  
\end{proof}

Now let $\pi\times U$ be a representation  of $A\times_\alpha G$ such
that $\pi$ is faithful. The representation $((\pi\times
U)\otimes\lambda)\circ\widehat\alpha$ has the form $(\pi\otimes
1)\times(U\otimes\lambda)$, and hence the induced
representation of $(A\times_\alpha G,G,\widehat\alpha)$ is given
by $\Ind(\pi\times U)=((\pi\otimes
1)\times(U\otimes\lambda),1\otimes M)$. It is trivial to check
that restricting $1\otimes M$ to $C_0(G/H)$
gives a covariant representation $((\pi\otimes
1)\times(U\otimes\lambda),1\otimes M|)$ of $(A\times_\alpha
G,G/H,\widehat\alpha)$, and hence we have a representation
$(\pi\otimes M|)\times(U\otimes\lambda)$ of the full crossed
product $C_0(G/H,A)\times_{\alpha\otimes\tau}G$ on
$\H_\pi\otimes L^2(G)$. Because we know from the Lemma that 
$\Ind(\pi\times U)$ is
faithful on $(A\times_\alpha G)\times_{\widehat\alpha} G$, the image of
$C_0(G/H,A)\times_{\alpha\otimes\tau}G$ is precisely (one realisation
of) the reduced crossed product
$(A\times_\alpha G)\times_{\widehat\alpha,r}G/H$. We shall call this
representation the \emph{regular representation of
$C_0(G/H,A)\times_{\alpha\otimes\tau}G$ induced from $(\pi,U)$.} As we
shall see, this representation is not always faithful.

\begin{prop}\label{reg-rep}
Suppose $(\pi,U)$ is a covariant representation of $(A,G,\alpha)$ on
$\H$ and $\pi$ is faithful. Then the regular representation $(\pi\otimes
M|)\times(U\otimes\lambda)$ induces an isomorphism of
$(A\otimes C_0(G/H))\times_{\alpha\otimes\tau,r}G$ onto $(A\times_\alpha
G)\times_{\widehat\alpha,r}G/H$.
\end{prop}

\begin{proof}
In view of the preceding remarks, it is enough to prove that the kernel
of $(\pi\otimes
M|)\times(U\otimes\lambda)$ is precisely the kernel of a regular
representation of $(A\otimes C_0(G/H))\times_{\alpha\otimes\tau}G$. We
know from Lemma \ref{ind-faith} that $\Ind(\pi\times U)$ is faithful on
$(A\times_\alpha G)\times_{\widehat\alpha}G$, hence 
 Lemma \ref{fund-lem} implies that $(\pi\otimes M,U\otimes\lambda)$ is
faithful on $(A\otimes C_0(G))\times_{\alpha\otimes\tau}G$.

The inclusion of $C_0(G/H)$ in $M(C_0(G))$ induces a homomorphism
$\phi$ of the crossed product $(A\otimes
C_0(G/H))\times_{\alpha\otimes\tau}G$ into
$M((A\otimes C_0(G))\times_{\alpha\otimes \tau}G)$. The regular
representation $\Ind(\pi\otimes M)$ is faithful on $(A\otimes
C_0(G))\times_{\alpha\otimes\tau}G$, and the composition
$\Ind(\pi\otimes M)\circ\phi$ is the regular
representation induced by the faithful representation $\pi\otimes M|$ of
$A\otimes C_0(G/H)$. Thus the kernel of $\phi$ is the kernel of the
regular representation, and $\phi$ induces an injection of $(A\otimes
C_0(G/H))\times_{\alpha\otimes\tau,r}G$ into $M((A\otimes
C_0(G))\times_{\alpha\otimes \tau}G)$. Composing this injection with the
faithful representation $(\pi\otimes M)\times (U\otimes \lambda)$ gives
a faithful representation of $(A\otimes
C_0(G/H))\times_{\alpha\otimes\tau,r}G$; but 
$((\pi\otimes M)\times
(U\otimes \lambda))\circ\phi=(\pi\otimes M|)\times (U\otimes \lambda)$,
so the result follows.
\end{proof}

\begin{cor}\label{reg-rep-faith}
We have $A\times_\alpha H=A\times_{\alpha,r}H$ if and only if, whenever
$(\pi,U)$ is a covariant representation of $(A,G,\alpha)$ with $\pi$
faithful, the regular representation of $(A\otimes
C_0(G/H))\times_{\alpha\otimes\tau}G$ induced from $(\pi,U)$ is faithful.
\end{cor}

\begin{proof}
Recall from \cite{QS-RH} that $A\times_\alpha H=A\times_{\alpha,r}H$ if
and only if
the imprimitivity algebra $(A\otimes
C_0(G/H))\times_{\alpha\otimes\tau}G$ is isomorphic to $(A\otimes
C_0(G/H))\times_{\alpha\otimes\tau,r}G$.
\end{proof}

\begin{rmk}\label{not-faith}
Applying this result with $H$ normal and amenable gives \cite[Lemma
3.2]{KQ-IC}, albeit only for dual coactions
(cf. also \cite[Proposition 7]{ManJF}). Taking $H=G$, $A=\C$
and $G$ nonamenable shows that the regular representation in 
\propref{reg-rep}\ is not always faithful.
\end{rmk}

Restricting the action on the left of the bimodule of
Proposition \ref{R-bimod-prop} gives a right-Hilbert
$C_0(G,A)\times G-C_0(G/H,A)\times G$ bimodule, which by
Lemma \ref{fund-lem} we can view as a right-Hilbert
$A\times G\times G-C_0(G/H,A)\times G$ bimodule $Z_{G/H}^G(A\times G)$.
Using this, we can induce a covariant representation
$(\pi\times U,\mu)$ of $(A\times_\alpha
G,G/H,\widehat\alpha)$ to a representation $\Ind_{G/H}^G(\pi\times
U,\mu)$ of $(A\times_\alpha G)\times_{\widehat\alpha}G$, acting in a
completion of $Z_{G/H}^G\otimes\H_\pi$. Since the isomorphism
of $(A\times_\alpha G)\times_{\widehat\alpha}G$ with
$C_0(G,A)\times_{\alpha\otimes\tau}G$ carries the double dual
action into the action of $G$ used in \S1, we deduce from
Proposition \ref{R-bimod-prop} the following
representation-theoretic imprimitivity theorem:

\begin{prop}\label{imp-homo-sp}
Suppose $\alpha:G\to\Aut A$ is an action of a locally compact
group on a $C^*$-algebra $A$ and $H$ is a closed subgroup of
$G$. A representation $(\rho\times V)\times \nu$ of
$(A\times_\alpha G)\times_{\widehat\alpha}G$ is induced from a
covariant representation of $(A\times_\alpha
G,G/H,\widehat\alpha)$ if and only if there is a representation
$U$ of $H$ on $\H_\rho$ such that $((\rho\times V)\times
\nu,U)$ is a covariant representation of $((A\times_\alpha
G)\times_{\widehat\alpha}G, H,\widehat{\widehat\alpha}|)$. 
{\rm(}That is, if and only if 
the range of $U$ commutes with the ranges of $V$ and $\nu$,
and $\nu(\sigma_s(f))=U_s\nu(f)U_s^*$ for $s\in H$, $f\in
C_0(G)$.{\rm)}
\end{prop}

\begin{rmk}
From \cite{QS-RH}, we know that the imprimitivity bimodule of Proposition
\ref{R-bimod-prop} has as a (possibly proper) quotient a
$C_0(G,A)\times_r(G\times H)$--$C_0(G/H,A)\times_r G$ imprimitivity
bimodule $Z_r$. Since
\begin{eqnarray*}
C_0(G,A)\times_r(G\times H)&\cong&(C_0(G,A)\times_r G)\times_r H\\
&\cong&(C_0(G,A)\times G)\times_r H\\
&\cong&(A\times_\alpha G\times_{\widehat\alpha}G)\times_r H,
\end{eqnarray*}
we can by Proposition \ref{reg-rep}  realise $Z_r$ as a
right-Hilbert $(A\times_\alpha G)\times_{\widehat\alpha}G$--$(A\times_\alpha
G)\times_{\widehat\alpha,r}G/H$ bimodule, and use it to induce
representations from the reduced crossed product. We shall see in
Theorem \ref{DK-thm} that this induction process agrees with the one
studied in \cite{ManJF,KQ-IC} for normal $H$.
\end{rmk}

\section{Induction and duality}

In this section we show that, modulo duality, our induction process for
dual systems is the inverse of Green induction. Before stating our
theorem, we describe the three bimodules involved.

Consider an action $\alpha:G\to \Aut A$ and 
a closed, not-necessarily-normal subgroup $H$ of $G$. 
Recall from \cite{GreLS}\ that 
 $C_c(G,A)$ 
can be completed to a $\CGmHAxG$ -- $\AxaH$ \ib\ $X_H^G(A)$. 
We use the pre-\ib\ structure on $C_c(G,A)$ given for 
$f\in C_c(G\x G/H,A)$, $x$ and $y$ in $C_c(G,A)$, and
$g\in C_c(H,A)$ as follows:
\begin{eqnarray}
f\cdot x(r) &=& \int_G f(s,rH) \alpha_s(x(\inv{s}r))\, 
\Delta_G(s)^{\frac12}\,
ds\label{gre-ib-a}\\
x\cdot g(r) &=& \int_H x(rt)\alpha_{rt}(g(\inv{t}))\,
\Delta_H(t)^{-{\frac12}}\, 
dt\label{gre-ib-b}\\
\lip{C_0(G/H,A)\x G}{x}{y}(s,rH) &=& \int_H
x(rt)\alpha_s(y(\inv{s}rt)^*)\, \Delta_G(s)^{-{\frac12}}\,
dt\label{gre-ib-c}\\
\rip{A\x H}{x}{y}(t) &=& \int_G \alpha_s(x(\inv{s})^*y(\inv{s}t))\,
\Delta_H(t)^{-{\frac12}}\, 
ds.\label{gre-ib-d}
\end{eqnarray}
These actions and inner products, and in particular the modular
functions, come straight
from the symmetric imprimitivity theorem (see \S1),
with $K=G$ and $L=H$ acting on $P=G$ by left and right multiplication,
$\sigma=\alpha$, and $\rho=\id$.  

Recall from \S 1 that, in the case $H=\{e\}$, the action $\ahh$ of $G$ on
Green's imprimitivity algebra $C_0(G,A)\times G$ 
is given for $f\in C_c(G\x G,A)$ by
$
\ahh_t(f)(r,s) = f(r,st)
$. The imprimitivity bimodule $X_{\{e\}}^G(A)$ also
admits an action
$\gamma$ of
$G$, given for $x\in C_c(G,A)$ by
$\gamma_t(x)(s) = x(st)$,
and by \cite[Theorem~1]{EchME}\ this gives an equivariant Morita
equivalence $(X_{\{e\}}^G(A),\gamma) $between $(C_0(G,A)\times G,G,\ahh)$
and
$(A,G,\alpha)$.   Thus for any closed subgroup $H$ of $G$ we have
a $(C_0(G,A)\times_{\alpha\otimes\tau} G)\times_{\ahh} H$--$A\times_{\alpha}
H$ imprimitivity bimodule
$X_{\{e\}}^G(A)\times H$, with dense submodule $C_c(H\times G,A)$
\cite{CMW-CP}. For $f\in C_c(H\x G\x G,A)$, $x$ and $y$ in
$C_c(H\x G,A)$, and $g\in C_c(H,A)$, the actions and
inner products are as follows:
\begin{eqnarray}
\label{XH-eqn-a}
f\cdot x(h,r) &=& \int_G\int_H f(k,u,r)
\alpha_u(x(\inv{k}h,\inv{u}rk))\, \Delta_G(u)^{\frac12}\, dk\, du\\
\label{XH-eqn-b}
x\cdot g(h,r) &=& \int_H x(k,r)
\alpha_{rk}(g(\inv{k}h))\, 
dk\\
\label{XH-eqn-c}
\lip{L}{x}{y}(h,r,s) &=& \int_H x(k,s)
\alpha_r(y(\inv{h}k,\inv{r}sh)^*)\,
\Delta_H(\inv{h}k)\Delta_G(r)^{-{\frac12}}\, dk\\
\label{XH-eqn-d}
\rip{R}{x}{y}(h) &=& \int_G\int_H
\alpha_s(x(\inv{k},\inv{s}k)^*y(\inv{k}h,\inv{s}k))\,
\Delta_H(k)^{-1}\, dk\, ds.
\end{eqnarray}

As in the previous section, we denote by $Z_{G/H}^G(A\times G)$
the bimodule of Proposition \ref{R-bimod-prop} viewed as an 
$(A\times_{\alpha} G\times_{\ah}G)\times_{\ahh}
H-C_0(G/H,A)\times_{\alpha\otimes\tau} G$ imprimitivity bimodule.

\begin{thm}\label{HGK-thm}
Let $\alpha:G\to\Aut A$ be an action of a locally compact group $G$ on a
$C^*$-algebra $A$, and let
$H$ be a closed  subgroup of $G$.  Then 
\[
Z_{G/H}^G(A\times G)\otimes_{C_0(G/H,A)\times G} X_H^G(A)\cong 
X_{\{e\}}^G(A)\times H
\]
as $(A\times_{\alpha} G\times_{\ah}G)\times_{\ahh}
H-C_0(G/H,A)\times_{\alpha\otimes\tau} G$ imprimitivity bimodules. 
\end{thm}

For the proof, we shall need the special case of the following lemma in
which $\psi_A$ and $\psi_B$ are the identity; the general case will be
used in \S4.
 
\pagebreak
\begin{lem}\label{suff-map-lem}
Suppose that $_AX_B$ and $_CY_D$ are imprimitivity bimodules, let
$\psi_A:A\to C$, $\psi_B:B\to D$ be surjective homomorphisms, and let 
$J=\ker\psi_A, I=\ker\phi_B$.
If $\psi_X:X\to Y$ is a linear map satisfying
\begin{align*}
\psi_X(a\cdot x)&=\psi_A(a)\cdot\psi_X(x)\\
\psi_X(x\cdot b)&=\psi_X(x)\cdot\psi_B(b)\\
\rip{D}{\psi_X(x)}{\psi_X(y)}&=\psi_B(\rip{B}{x}{y}).
\end{align*}
Then $\ker \psi_X=X\cdot I$ and 
$(\psi_A,\psi_X,\psi_B)$ factors through an imprimitivity bimodule
isomorphism of $_{A/J}(X/X\cdot I)_{B/I}$ onto $_CY_D$.
\end{lem}
\begin{proof}
We have
\begin{eqnarray*}
\lip{C}{\psi_X(x)}{\psi_X(y)}\cdot\psi_X(z) & = &
\psi_X(x)\cdot\rip{D}{\psi_X(y)}{\psi_X(z)}\\
 & = & \psi_X(x)\cdot\psi_B(\rip{B}{y}{z})\\
 & = & \psi_X(x\cdot\rip{B}{y}{z})\\
 & = & \psi_X(\lip{A}{x}{y}\cdot z)\\
 & = & \psi_A(\lip{A}{x}{y})\cdot\psi_X(z).
\end{eqnarray*}
Since $\psi_A$ and $\psi_B$ are surjective, it follows that
$\overline{\psi_X(X)}$ is a full $C$--$D$ submodule of $_CY_D$.
Thus, by \cite[Theorem~3.1]{RieUR}, $\psi_X(X)$ is dense in $Y$.
Then the above computations imply that
$(\psi_A,\psi_X,\psi_B)$ is an imprimitivity bimodule homomorphism
which factors through an injective imprimitivity bimodule homomorphism
$(\psi_{A/J}, \psi_{X/X\cdot I}, \psi_{A/I})$
of \mbox{$_{A/J}(X/X\cdot I)_{B/I}$} into $_CY_D$ by \cite[Lemma 2.7]{ER-ST}.
Since $\psi_{X/X\cdot I}$ is isometric, it follows that
$\psi_X(X)$ is complete. Hence $\psi_X(X)=Y$.
\end{proof}

\begin{proof}[Proof of Theorem \ref{HGK-thm}.]
We work with the dense subalgebras 
\[
C_c(H\times G\times G,A)\subseteq (A\times_{\alpha} G\times_{\ah}
G)\times_{\ahh} H 
\quad\text{and}\quad 
C_c(H,A)\subseteq A\times_{\alpha} H,
\]
and the
dense submodules
$$
C_c(G\times G,A)\subseteq Z_{G/H}^G(A\times G),
\quad C_c(G,A)\subseteq X_H^G(A),\
\text{and}\quad C_c(H\times G,A)\subseteq X_{\{e\}}^G(A)\times H.  
$$

Fix $(f,x)$ in $C_c(G\times G,A)\times C_c(G,A)$ and 
suppose $E_{f_1}$, $E_{f_2}$, and $E_x$ are compact sets
such that
$\supp(f)\subseteq E_{f_1}\x E_{f_2}$ and $\supp(x) \subseteq E_x$;
then the map 
$F_{f,x}\colon H\x G\x G\to A$ defined by
$$
F_{f,x}(h,s,t) = f(t,s) \alpha_t(x(t^{-1}sh)) \Delta_H(h)^{-\frac12}
\Delta_G(t)^{\frac12}
$$
is continuous and has support in 
$(E_{f_2}^{-1} E_{f_1} E_x)\cap H \times E_{f_2} \times E_{f_1}$.
It follows that the map $(h,s)\mapsto \int_G F_{f,x}(h,s,t)\, dt$ is in
$C_c(H\x G,A)$.  The pairing 
which sends $(f,x)$ to this element of 
$C_c(H\x G,A)$ is bilinear, 
and so we have a well-defined map
$\psi$ of $C_c(G\x G,A)\odot C_c(G,A)$ into $C_c(H\x G,A)$ given by
$$
\psi(f\ox x)(h,s) 
= \int_G f(t,s) \alpha_t(x(t^{-1}sh)) 
\Delta_H(h)^{-\frac12} \Delta_G(t)^{\frac12}\, dt.
$$

The following 
calculations verify that $\psi$ preserves both
actions and the right inner product.
For $g\in C_c(H\times G\times G,A)$ and $f\otimes x\in C_c(G\times
G,A)\odot C_c(G,A)$:
\begin{eqnarray*}
\psi(g\cdot f\otimes
x)(h,s)
 & = & \int_G g\cdot f(t,s) \alpha_t(x(t^{-1}sh))
\Delta_H(h)^{-\frac12}\Delta_G(t)^{\frac12}\, dt\\
 & \stackrel{t\mapsto ut}{=} & 
\int_G\int_G\int_H g(k,u,s) \alpha_u\bigl( f(t,u^{-1}sk)
\alpha_t(x(t^{-1}u^{-1}sh)) \\
 &   & \quad\quad\quad\quad\quad\quad\quad\quad\quad\quad
\Delta_H(k^{-1}h)^{-\frac12}
\Delta_G(t)^{\frac12}\bigr) \Delta_G(u)^{\frac12}\, dk\, du\, dt\\
 & = & \int_G\int_H g(k,u,s) \alpha_u(\psi(f\otimes
x)(k^{-1}h,u^{-1}sk)) \Delta_G(u)^{\frac12}\, dk\, du\\
& \stackrel{(\ref{XH-eqn-a})}{=} & 
g\cdot \psi(f\otimes x)(h,s).
\end{eqnarray*}
For $f\otimes x\in C_c(G\times G,A)\odot C_c(G,A)$ and $g\in
C_c(H,A)$:
\begin{eqnarray*}
\psi(f\otimes x\cdot g)(h,s)
 & = & \int_G f(t,s) \alpha_t(x\cdot g(t^{-1}sh))
\Delta_H(h)^{-\frac12}\Delta_G(t)^{\frac12}\, dt\\
 & = & \int_G\int_H f(t,s) \alpha_t(x(t^{-1}shk))
\alpha_{shk}(g(k^{-1})) \Delta_H(hk)^{-\frac12}\Delta_G(t)^{\frac12}\,
dk\, dt\\
 & \stackrel{k\mapsto h^{-1}k}{=} & 
\int_G\int_H f(t,s) \alpha_t(x(t^{-1}sk))
\alpha_{sk}(g(k^{-1}h)) \Delta_H(k)^{-\frac12}\Delta_G(t)^{\frac12}\,
dk\, dt\\
 & = & \int_H \psi(f\otimes x)(k,s) \alpha_{sk}(g(k^{-1}h))\, dk\\
 & \stackrel{(\ref{XH-eqn-b})}{=} & 
\psi(f\otimes x)\cdot g\, (h,s).
\end{eqnarray*}
For $f\otimes x$ and $g\otimes y$ in $C_c(G\times G,A)\odot C_c(G,A)$:
\begin{align*}
\langle f\otimes &x, g\otimes y\rangle_{A\times H}(h)
=  \langle x, \langle f,g\rangle_{C_0(G/H,A)\times G}\cdot
y\rangle_{A\times H}(h) \\
 & \stackrel{(\ref{gre-ib-d})}{=} 
\int_G \alpha_s( x(s^{-1})^* \langle f,g\rangle_R\cdot y(s^{-1}h))
\Delta_H(h)^{-\frac12}\, ds\\
 & \stackrel{(\ref{gre-ib-a})}{=} 
\int_G\int_G \alpha_s\bigl( x(s^{-1})^* \langle f,g\rangle_R(t,s^{-1}H) 
\alpha_t(y(t^{-1}s^{-1}h)) \Delta_G(t)^{\frac12}\bigr)
\Delta_H(h)^{-\frac12}\, dt\, ds\\
 & \stackrel{(\ref{K-eqn-d})}{=} 
\int_G\int_G\int_G\int_H \alpha_s\bigl( x(s^{-1})^* \alpha_u(
f(u^{-1},u^{-1}s^{-1}k)^* g(u^{-1}t,u^{-1}s^{-1}k)) \\
 &    \quad\quad\qquad\quad\quad\qquad\quad\quad\qquad
\Delta_G(u^{-1})
\alpha_t(y(t^{-1}s^{-1}h)) \Delta_G(t)^{\frac12}\bigr)
\Delta_H(h)^{-\frac12}\, dk\, du\, dt\, ds\\
 & \stackrel{u\mapsto u^{-1}}{=} 
\int_G\int_G\int_G\int_H \alpha_s(x(s^{-1})^*)
\alpha_{su^{-1}}(f(u,us^{-1}k)^*) \alpha_{su^{-1}}(g(ut,us^{-1}k))\\
 &    \quad\quad\qquad\quad\quad\qquad\quad\quad\qquad\quad\quad\qquad
\alpha_{st}(y(t^{-1}s^{-1}h))
\Delta_G(t)^{\frac12}\Delta_H(h)^{-\frac12}\, dk\, du\, dt\, ds\\
&\stackrel{t\mapsto u^{-1}t}{\stackrel{s\mapsto su}{=}} 
\int_G\int_G\int_G\int_H \alpha_s\bigl( (
f(u,s^{-1}k)\alpha_u(x(u^{-1}s^{-1}))
\Delta_G(u)^{\frac12} )^*g(t,s^{-1}k) \\
 &    \quad\quad\quad\quad\quad\quad\quad\quad\quad\quad\quad\quad
\quad\quad\quad\quad \alpha_t(y(t^{-1}s^{-1}h))
\Delta_H(h^{-1})^{\frac12}\Delta_G(t)^{\frac12} \bigr)\, 
dk\, du\, dt\, ds
\end{align*}
\pagebreak
\begin{align*}
\phantom{\langle f\otimes}
 &=  \int_G\int_H \alpha_s( \psi(f\otimes x)(k^{-1},s^{-1}k)^*
\psi(g\otimes y)(k^{-1}h,s^{-1}k)) 
\Delta_H(k^{-1})\, dk\, ds\\
 & \stackrel{(\ref{XH-eqn-d})}{=} 
\langle \psi(f\otimes x),\psi(g\otimes y)\rangle_{A\times H}(h).
\end{align*}
It follows that $\psi$ 
extends to a linear map of 
$Z_{G/H}^G(A\times G)\otimes_{C_0(G/H,A)\times G}X_H^G(A)$ into 
$X_{\{e\}}^G(A)\times H$
which also preserves the actions and
right inner product, and which therefore by \lemref{suff-map-lem}\ is
actually an isomorphism of the \ib s.  
\end{proof}

\begin{cor}\label{ind-res-cor}
Let $\alpha:G\to \Aut A$ be an action of a locally compact group on a
$C^*$-algebra, and let
$H$ be a closed subgroup of $G$.  Then we have a commutative diagram
\begin{equation*}
\begin{diagram}
\node{\Rep A\times_{\alpha} H}
	\arrow{s,l}{\Res_{\{e\}}^H}
	\arrow[2]{e}
\node[2]{\Rep C_0(G/H,A)\times_{\alpha\otimes\tau} G}
	\arrow{s,r}{\Ind_{G/H}^G}\\
\node{\Rep A}
	\arrow[2]{e}
\node[2]{\Rep \AxGxG}
\end{diagram}
\end{equation*}
in which the horizontal arrows are the bijections induced by the Green
bimodules $X_H^G(A)$ and $X_{\{e\}}^G(A)$.
\end{cor}

\begin{proof}
We shall show rather more: each arrow is implemented by a
right-Hilbert bimodule, so the two compositions are implemented by
tensor products of these bimodules, and we shall show that
\[
Z_{G/H}^G(A\times G)\otimes_{C_0(G/H,A)\times G} X_H^G(A)
\cong 
X_{\{e\}}^G(A)\otimes_A (A\times_{\alpha} H)
\]
as right-Hilbert $\AxGxG$ -- $A\times_{\alpha}H$ bimodules.  But the
bimodule on the right-hand side  is isomorphic 
to the right-Hilbert $\AxGxG$ -- $A\times_{\alpha}H$ bimodule
$X_{\{e\}}^G(A)\times H$ by a special case of 
\cite[Lemma~5.7]{KQR-DR}, 
so the isomorphism follows from \thmref{HGK-thm}.
\end{proof}


\section{Comparison with Mansfield's bimodule}\label{mansfield-sec}

Here we compare our inducing process for dual coactions
with that of \cite{KQ-IC}, which extends Mansfield's process to nonamenable
subgroups. For each coaction $\delta:B\to M(B\otimes C^*(G))$ 
and normal subgroup $N$, \cite[Theorem 3.3]{KQ-IC} provides an imprimitivity
bimodule
$Y_{G/N}^G$ between the reduced crossed products
$(B\times_{\delta}G)\times_{\widehat{\delta},r}N$ and $B\times_{\delta,r}G/N$.

We consider an action $\alpha:G\to \Aut A$, the dual coaction
$\widehat\alpha$ on $A\times_\alpha G$, and a closed normal subgroup $N$ of
$G$. To define the reduced crossed product, we fix a faithful representation
$\pi$ of $A$ on $\H$, and use the covariant representation
$\Ind\pi:=\widetilde\pi\times(1\otimes\lambda)$ of $A\times_\alpha G$ on
$\H\otimes L^2(G)$. 
Note that $\widetilde\pi$ is faithful, so Proposition \ref{reg-rep}
gives us a faithful representation $(\widetilde\pi\otimes M|)\times
(1\otimes\lambda\otimes\lambda)$ of  $(A\otimes
C_0(G/N))\times_{\alpha\otimes\tau,r}G$ onto $(A\times_\alpha
G)\times_{\widehat\alpha,r}G/N$.

We now recall the construction of the bimodule from \cite{ManJF,KQ-IC}. Consider
the map $\varphi\colon C_c(G)\to C_c(G/N)$ defined by 
$$
\varphi(f)(rN) = \int_N f(rh)\, dh.
$$
Then $\D_N$ is a *-subalgebra of
$B(\H\ox\LtwoG\ox\LtwoG)$ containing in particular the elements of the
form
$$\Ind\pi\ox\lambda(\ah(\ah_u(b)))(1\ox1\ox M(\varphi(f)))
\quad\text{and}\quad
(1\ox1\ox M(\varphi(f)))\Ind\pi\ox\lambda(\ah(\ah_u(b))),$$
for $b\in \AxaG$, $u\in A_c(G)$, and $f\in C_c(G)$.  
(By definition, $\ah_u$, is the composition of $\ah$ with the slice map
$S_u:=\id\otimes u:M(\AxaG\ox C^*(G))\to M(\AxaG)$.)
$\D$~is by definition $\D_{\{e\}}$.
Mansfield shows that there is a well-defined map $\varPsi \colon\D\to\D_N$
such that
\begin{equation}\label{varPsi-def-eqn-a}
\varPsi \left(\Ind\pi\otimes\lambda(\ah(\ah_u(b)))\right) 
= \Ind\pi\otimes\lambda(\ah(\ah_u(b))),\ \text{ and}
\end{equation}
\begin{equation}\label{varPsi-def-eqn-b}
\varPsi \left(1\otimes1\otimes M(f)\right) 
= 1\otimes 1\otimes M(\varphi(f)).
\end{equation}
Then $\D$ has a $\D_N$-valued pre-inner product given by
$$
\rip{\D_N}{x}{y} = \varPsi (x^*y).
$$
With left action of $\I_N\subseteq C_c(N,\D)$ 
given by
\begin{equation}\label{left-D-act-eqn}
f\cdot x = \int_N f(n) \ahh_n(x)\, \Delta_N(n)^{\frac12}\, dn
\end{equation}
and right action of $\D_N$ given by
$
x\cdot z = xz
$,
$\D$ becomes an $\I_N$--$\D_N$ pre-\ib, whose completion
$Y_{G/N}^G(A\times G)$ is an 
$(A\times_{\alpha}G\times_{\widehat{\alpha}}G)\times_{\ahh,r}N-
(A\times_{\alpha}G)\times_{\ah,r}G/N$ \ib
 (\cite[Theorem 3.3]{KQ-IC}, \cite[Theorem~27]{ManJF}).
Recall that our bimodule $Z_{G/N}^G(A\times_{\alpha}G)$ is
an imprimitivity bimodule between the full crossed products
$(A\times_{\alpha}G\times_{\widehat{\alpha}}G)\times_{\ahh}N$
and $(A\otimes C_0(G/N))\times_{\alpha\otimes\tau}G$.

\begin{thm}\label{DK-thm}
Let $\alpha:G\to\Aut A$ be an action of a locally compact group $G$ on a
$C^*$-algebra $A$, and let
$N$ be a closed normal subgroup of $G$. Let
\[
\Upsilon:=(\widetilde\pi\times M|)\times(1\otimes\lambda\otimes\lambda):
(A\otimes C_0(G))\times_{\alpha\otimes\tau} G
\to (A\times_{\alpha} G)\times_{\ah,r} G/N,
\]
and let $\Phi:(A\times_\alpha
G\times_{\widehat\alpha}G)\times_{\widehat{\widehat\alpha}}N\to
\L(Y_{G/N}^G)$ be the extension of the left action (\ref{left-D-act-eqn}).
Then there exists a linear map
$\Theta$ of $Z_{G/N}^G(A\times G)$ onto $Y_{G/N}^G(A\times G)$ such that
$(\Phi,\Theta,\Upsilon)$ is a surjective imprimitivity bimodule
homomorphism. In particular, if 
$I=\ker\Upsilon$, then 
$$Z_{G/N,r}^G(A\times G):=Z_{G/N}^G(A\times G)/\big(Z_{G/N}^G(A\times G)\cdot
I\big)
\cong Y_{G/N}^G(A\times G)$$
as
$(A\times_{\alpha}G\times_{\ah}G)\times_{\ahh,r}N$--$(A\times_{\alpha}
G)\times_{\ah,r} G/N$ imprimitivity bimodules.
\end{thm}

\begin{proof}
It is sufficient to produce a linear
map
$\Theta$ of a dense subspace $Z_0\subset Z$ into $\D$ 
such that
\begin{equation}\label{DK-eqn-a}
\Theta(f\cdot x) = \Phi(f)\cdot\Theta(x),
\end{equation}
\begin{equation}\label{DK-eqn-b}
\Theta(x\cdot g) = \Theta(x)\cdot\Upsilon(g),\ \text{ and}
\end{equation}
\begin{equation}\label{DK-eqn-c}
\rip{(A\times_{\alpha}G)\times_{\ah,r}G/N}{\Theta(x)}{\Theta(y)} 
  = \Upsilon\big(\rip{C_0(G/N,A)\times_{\alpha\otimes\tau}G}{x}{y}\big),
\end{equation}
for $f\in C_c(N\x G\x G,A)$, $g\in C_c(G\x
G/N,A)$, and $x,y\in Z_0$;
for then $\Theta$ extends to a linear map of
$Z_{G/N}^G(A\times G)$ into $Y_{G/N}^G(A\times G)$ which also satisfies
\eqref{DK-eqn-a}--\eqref{DK-eqn-c}, and hence
factors through an imprimitivity
bimodule isomorphism of $Z_{G/N,r}^G(A\times G)$
onto  $Y_{G/N}^G(A\times G)$ by Lemma \ref{suff-map-lem}.

Let  $\Theta$ be the restriction of
$(\tilde{\pi}\times M)\times (1\otimes \lambda\otimes \lambda)$
to $C_c(G\times G,A)\subseteq (A\otimes C_0(G))\times_{\alpha\otimes
\tau}G$, and let
$$
Z_0 = \spn\{ a\otimes z\otimes f \mid a\in A;\ z,f\in C_c(G) \} \\
\subseteq C_c(G\times G,A).
$$
By Lemma \ref{fund-lem} we have
\begin{equation}\label{elem-theta-eqn}
\Theta(a\ox z\ox f) 
= (1\ox1\ox M_f)\Ind\pi\ox\lambda(\ah(a\ox z)).
\end{equation}
Choosing $u\in A_c(G)$ to be identically $1$ on $\supp(z)$, we have
$$
\ah(a\ox z) = \ah(a\ox uz) = \ah(\ah_u(a\ox z)),
$$
because $S_u(\ah(g))$ is the pointwise product $ug$
\cite[Lemma~1.3]{RaeCR}. 
Thus $\Theta$ maps $Z_0$ into $\D$. 
To see that $Z_0$ is dense in $Z_{G/N}^G(A\times G)$, note
that the inductive limit topology dominates the
\ib\ topology on $C_c(G\times G,A)$ \cite[p.~374]{RaeIC}.

To verify \eqnref{DK-eqn-a}, notice that for $f\in
C_c(N\times G\times G, A)$
and $x\in Z_0$, 
\eqnref{K-eqn-a}\ may be re-written as
$$
f\cdot x = \int_N f(n)*\ahh_n(x)\, \Delta_N(n)^{\frac12}\, dn,
$$
where $*$ denotes the 
multiplication in the subalgebra $C_c(G\x G,A)$ of $\AxGxG$.
Since $\Theta$ is multiplicative with respect to this operation,
this gives
\begin{eqnarray*}
\Theta(f\cdot x) 
 & = & \Theta\left( \int_N f(n)*\ahh_n(x)\, \Delta_N(n)^{\frac12}\, dn\right)\\
 & = & \int_N \Theta(f(n)) \Theta(\ahh_n(x))\, \Delta_N(n)^{\frac12}\, dn\\
 & = & \int_N \Phi(f)(n) 
\ahh_n(\Theta(x))\, \Delta_N(n)^{\frac12}\, dn\\
 & = & \Phi(f)\cdot\Theta(x).  
\end{eqnarray*}

To verify \eqnref{DK-eqn-b} and \eqnref{DK-eqn-c}, 
we first let $a\otimes z\otimes f\in Z_0$ and 
$\xi\in L^2(G\times G,\H)\cong \H\otimes L^2(G)\otimes
L^2(G)$, and compute: 
\begin{eqnarray*}
(\Theta(a\otimes z\otimes f)\xi)(r,s) &
 \stackrel{(\ref{elem-theta-eqn})}{=}& 
\big((1\otimes 1\otimes M_f)\Ind\pi\otimes\lambda(\ah(a\otimes
z))\xi\big)(s,t)\\  
& {=} & 
\big((1\otimes 1\otimes M_f)(\tilde{\pi}\otimes 1)(a)(1\otimes
\lambda\otimes
\lambda)(z)\xi\big)(s,t)\\
&=& \int_G \pi\big(\alpha_{r^{-1}}(az(t) f(s))\big)\xi(t^{-1}r,t^{-1}s)\, dt.
\end{eqnarray*}
Thus, since $Z_0$ is inductive-limit dense in $C_c(G\times G,A)$,
it follows that
for all $g$ in the subalgebra  $C_c(G\times G,A)$ of 
$(A\otimes C_0(G))\times_{\alpha\otimes \tau}G$, we have
\begin{equation}\label{tilde-Theta-eqn}
 \Theta(g)\xi(s,t)=
\int_G\pi\big(\alpha_{r^{-1}}(g(t,s))\big)\xi(t^{-1}r,t^{-1}s)\,dt.
\end{equation}
Since $\Upsilon$ is the restriction of $\Theta$ to the image
of
$(A\otimes C_0(G/N))\times_{\alpha\otimes \tau}G$ in $M((A\otimes
C_0(G))\times_{\alpha\otimes\tau}G)$, where for the moment we
identify $\Theta$ with its extension to $(A\otimes
C_0(G))\times_{\alpha\otimes\tau}G$, it also follows that
\begin{equation}\label{tilde-Psi-eqn}
\Upsilon(g)\xi(s,t)=
\int_G\pi\big(\alpha_{r^{-1}}(g(t,sN))\big)\xi(t^{-1}r,t^{-1}s)\,dt.
\end{equation}
Notice that 
for $x\in Z_0$ and $g\in C_c(G\times G/N,A)$, 
Equation~\eqref{K-eqn-b}\ can be re-written as
$ x\cdot g = x*g$,
where $x*g$ denotes convolution of $x\in Z_0$ with $g\in C_c(G\times
G/N,A)$. Thus \eqref{DK-eqn-b} follows from
\[
\Theta(x\cdot g) 
 =\Theta(x*g)=\Theta(x)\Theta(g)=\Theta(x)\cdot\Upsilon(g).
\]

Before checking \eqnref{DK-eqn-c}, we need to do some background
calculations. First,  since
$\Theta$ is involutive on
$C_c(G\x G,A)$, we have for $x$ and $y$ in $Z_0$
\begin{equation*}
\rip{(A\times_{\alpha}G)\times_{\widehat{\alpha},r}G/N}{\Theta(x)}{\Theta(y)} =
\varPsi (\Theta(x)^*\Theta(y)) =
\varPsi (\Theta(x^**y));
\end{equation*}
thus to establish \eqnref{DK-eqn-c}, it is enough to verify
that 
$$
\varPsi (\Theta(x^**y)) = \Upsilon\left(\rip{\CGmNAxG}{x}{y}\right).
$$
Next, note that by \eqref{elem-theta-eqn}, 
\eqref{varPsi-def-eqn-a}, and \eqref{varPsi-def-eqn-b},
for $a\otimes z\otimes f\in Z_0$ we have
\[
\varPsi (\Theta(a\otimes z\otimes f)) = \Theta(a\otimes z\otimes
\varphi(f));
\]
thus we can compute:
\begin{eqnarray*}
\big(\varPsi (\Theta(a\otimes z\otimes f))\xi\big)(r,s)
&= &\big(\Theta(a\otimes z\otimes \varphi(f))\xi\big)(r,s)\\
& \stackrel{(\ref{tilde-Theta-eqn})}{=} & 
\int_G \pi(\alpha_{r^{-1}}(az(t)\varphi(f)(sN)))\xi(t^{-1}r,t^{-1}s)\, dt\\
& = & \int_G\int_N \pi(\alpha_{r^{-1}}(az(t)f(sn)))\xi(t^{-1}r,t^{-1}s)\, dn\, dt\\
& = & \int_G\int_N \pi(\alpha_{r^{-1}}(a\otimes z\otimes 
f(t,sn)))\xi(t^{-1}r,t^{-1}s)\, dn\, dt. 
\end{eqnarray*}
Hence by continuity we have
\begin{equation}\label{Psi-m-Theta-eqn}
\varPsi (\Theta(g)) = \int_G\int_N
\pi(\alpha_{r^{-1}}(g(t,sn)))\xi(t^{-1}r,t^{-1}s)\, dn\, dt
\end{equation}
for $g\in C_c(G\times G,A)$. 

Now to check \eqnref{DK-eqn-c}, we fix $x,y\in Z_0$  and compute:
\begin{eqnarray*}
\lefteqn{\big(\Upsilon\big(\rip{\CGmNAxG}{x}{y}\big)\xi\big)(r,s){\mbox{\qquad}}}\\
 & \stackrel{(\ref{tilde-Psi-eqn})}{=} & 
\int_G \pi\big(\alpha_{\inv{r}}(\rip{\CGmNAxG}{x}{y}(t,sN))\big)
\xi(\inv{t}r,\inv{t}s)\, dt\\
 & \stackrel{(\ref{K-eqn-d})}{=} & \int_G\int_G\int_N \pi\left(
\alpha_{\inv{r}u}(x(\inv{u},\inv{u}sn)^* y(\inv{u}t,\inv{u}sn))\right)
\xi(\inv{t}r,\inv{t}s)\Delta_G(u^{-1})\,
dn\, du\, dt\\
 & = & \int_G\int_N \pi\left( \alpha_{\inv{r}}\left( \int_G x^*(u,sn)
\alpha_u(y(\inv{u}t,\inv{u}sn))\,
du\right)\right)\xi(\inv{t}r,\inv{t}s)\, dn\, dt\\
 & = & \int_G\int_N \pi\big( \alpha_{\inv{r}}(x^**y(t,sn))\big)
\xi(\inv{t}r,\inv{t}s)\, dn\, dt\\
 & \stackrel{(\ref{Psi-m-Theta-eqn})}{=} & 
\big(\varPsi (\Theta(x^**y))\big)\xi(r,s).
\end{eqnarray*}
This completes the proof.
\end{proof}


\section{Appendix}

We prove the following weak version of
Mansfield's imprimitivity theorem  for the reduced crossed
product $B\times_{\delta,r}G/H$ of \S 2:

\begin{thm}\label{thm-genimp}
Let $\delta:B\to M(B\otimes C^*(G))$ be a nondegenerate coaction of
$G$ on $B$ and $H$ a closed subgroup of $G$. 
Then the reduced crossed product
$B\times_{\delta,r}G/H$ is  Morita equivalent to
$(B\times_{\delta}G)\times_{\widehat{\delta},r}H$.
\end{thm}

We saw at the end of \S2 that the Theorem is true for dual coactions,  
so we use the Morita equivalence of $\delta$ and
${\delta}\,\widehat{\widehat{\ }}$ to reduce to this case. There is one subtlety 
involved: if $\delta:B\to M(B\otimes C^*(G))$ is an arbitrary full coaction, 
it may not be true
that $\delta$ is Morita equivalent to ${\delta}\,\widehat{\widehat{\ }}$.
However, from Katayama's Duality Theorem \cite{KatTD}
we can deduce that this is true for nondegenerate reduced coactions (see Proposition
\ref{prop-Mor} below). 

We recall the definition of the reduction of a coaction
$\delta:B\to M(B\otimes C^*(G))$  from \cite{RaeCR,QuiFR}. Let 
$p:B\to B^r:=B/\ker j_B$ denote the quotient map. Then there is a well-defined
homomorphism 
$\delta^r:B^r\to M(B^r\otimes C_r^*(G))$ such that
$\delta^r\circ p=(p\otimes\lambda)\circ \delta$, and $\delta^r$ is a reduced
coaction of $G$ on $B^r$ which is nondegenerate if $\delta$ is (\cite[Lemma
3.1]{RaeCR}, \cite[Corollary 3.4]{QuiFR}). The canonical map $j_B$ factors through an
embedding $j_{B^r}$ of $B^r$ in $M(B\times_\delta G)$, and then $(B\times_\delta
G,j_{B^r},j_{C(G)})$ is a crossed product for the reduced system $(B^r,G,\delta^r)$.
Thus both reduced crossed products in Theorem
\ref{thm-genimp} depend only on the reduced system, and  Theorem
\ref{thm-genimp} will be a corollary of:

\begin{thm}\label{thm-genimp1}
Let $\delta:B\to M(B\otimes C_r^*(G))$ be a nondegenerate reduced 
coaction of $G$ on $B$ and assume that $B$ is represented 
faithfully and nondegenerately on a Hilbert space $\H$. Then 
$$B\times_{\delta,r}G/H=\clsp\{\delta(b)(1\otimes M_f): b\in B, f\in C_0(G/H)\}$$
is Morita equivalent to $(B\times_{\delta}G)\times_{\widehat{\delta},r}H$.
\end{thm}

From now on, all coactions will be reduced.
Recall that a Morita equivalence $(X,\dx)$ between two
cosystems $(A,G,\da)$ and $(B, G,\db)$  consists of an
$A-B$ imprimitivity bimodule $X$ together with a linear map
$\dx:X\to M(_{A\otimes C_r^*(G)}(X\otimes C_r^*(G))_{B\otimes C_r^*(G)})$ 
such that
$(\da, \dx,\db)$ is an  imprimitivity bimodule homomorphism, and such that
$\dx$ satisfies the coaction identity $(\dx\otimes \id_G)\circ \dx=
(\id_X\otimes \dg)\circ \dx$ (see \cite{ER-MI} for more details).

\begin{ex}\label{ex-stablyext}
(1) \emph{Stabilised coactions.}  Suppose that 
$\delta:B\to M(B\otimes C_r^*(G))$ is
a coaction. Let $\sigma:C_r^*(G)\otimes \K(\H)\to \K(\H)\otimes C_r^*(G)$
denote the flip map.
Then $\delta^s=(\id_B\otimes \sigma)\circ (\delta\otimes \id_{\K})$
is a coaction of $G$ on $B\otimes \K(\H)$, called the {\em stabilised
coaction} of $\delta$. 
Let $X:=B\otimes   \H$ viewed as an $B\otimes \K(\H)-B$ imprimitivity
bimodule.
Then the map
$\delta_X:= (\id_B\otimes \sigma_{\H})\circ (\delta\otimes \id_{\H})$
of $X$ into $M(X\otimes C_r^*(G))$
is a Morita equivalence for $\delta^s$ and $\delta$, where now
$\sigma_{\H}$ denotes the flip map between the imprimitivity bimodules
$_{C_r^*(G)\otimes \K(\H)}(C_r^*(G)\otimes \H)_{C_r^*(G)}$ and
$_{\K(\H)\otimes C_r^*(G)}(\H\otimes C_r^*(G))_{C_r^*(G)}$.\\
(2) \emph{Exterior equivalent coactions.} 
A $\delta$-one cocycle for a coaction $\delta:B\to M(B\otimes C_r^*(G))$
is a unitary $V\in UM(B\otimes
C_r^*(G))$ satisfying
$(\id_B\otimes \delta_G)(V)=(V\otimes 1)\big((\delta\otimes \id_G)(V)\big)$
and $V\delta(b)V^*(1\otimes z)\in B\otimes C_r^*(G)$ for all $b\in B, z\in
C_r^*(G)$ (see \cite[Definition 2.7]{LPRS-RC}). 
Then $\varepsilon=\Ad V\circ \delta$ is a coaction 
of $G$ on $B$. If $X=B$ is the trivial $B-B$ imprimitivity bimodule,
then $\delta_X:b\mapsto V\delta(b)$
is a Morita equivalence between $\varepsilon$ and $\delta$.
\end{ex}

\begin{prop}\label{prop-Mor}
If $\delta:B\to M(B\otimes C_r^*(G))$ is a nondegenerate reduced
coaction, then $\delta$ is Morita equivalent to the double
dual coaction ${\delta}\,\widehat{\widehat{\ }}$ of $G$ on
$(B\times_{\delta}G)\times_{\widehat{\delta},r}G$.
\end{prop}

\begin{proof}
It follows from \cite[Theorem~8]{KatTD} that there is an isomorphism of
$(B\times_{\delta}G)\times_{\widehat{\delta},r}G$ onto
$ B\otimes \K(L^2(G))$ carrying
${\delta}\,\widehat{\widehat{\ }}$ to the coaction 
$\Ad V\circ \delta^s$, where
$V=1\otimes W_G^*\in UM(B\otimes \K(L^2(G))\otimes C_r^*(G))$ is
a $\delta^s$-one cocycle.
Thus $\delta$ is Morita equivalent to
${\delta}\,\widehat{\widehat{\ }}$ by Example~\ref{ex-stablyext}.
\end{proof}

\begin{prop}\label{prop-meq}
If $(X,\dx)$ is a Morita equivalence for the cosystems
$(A, G,\da)$ and $(B,G,\db)$, and $H$ is a closed subgroup of $G$, then there
is an $A\times_{\da,r}G/H-B\times_{\db,r}G/H$ imprimitivity
bimodule $X\times_{\dx,r}G/H$.
\end{prop}
\begin{proof}
Let $L=\left(\begin{smallmatrix} A&X\\ \widetilde{X}&B\end{smallmatrix}\right)$
denote the linking algebra for $_AX_B$, and let 
$\dl=
\left(\begin{smallmatrix} \da&\dx\\
\delta_{\widetilde{X}}&\db\end{smallmatrix}\right)$ denote the 
corresponding coaction
of $G$ on $L$ (see \cite[Appendix]{ER-MI}).  We can represent $L$ faithfully on
$\H\oplus\K$ in such a way
that  the corners $A=pLp$ and $B=qLq$,
$p=\left(\smallmatrix 1&0\\ 0&0\endsmallmatrix\right)$,
$q=\left(\smallmatrix 0&0\\ 0&1\endsmallmatrix\right)$,
act faithfully and nondegenerately on
$\H$ and $\K$. Then 
$$L\times_{\dl,r}G/H=\clsp\{\dl(l)(1\otimes M_f):l\in L, f\in C_0(G/H)\},$$
and if $p\otimes 1,q\otimes 1$ denote the
projections of
$(\H\oplus\K)\otimes L^2(G)\cong(\H\otimes L^2(G))\oplus 
(\K\otimes L^2(G))$ onto its factors, then
$$(p\otimes 1)\big(L\times_{\dl,r}G/H\big)(p\otimes 1)= A\times_{\da,r}G/H,\ \text{
 and}$$
$$(q\otimes 1)\big(L\times_{\dl,r}G/H\big)(q\otimes 1)= B\times_{\db,r}G/H.$$
We claim that 
$$X\times_{\dx,r}G/H:=(p\otimes 1)\big(L\times_{\dl,r}G/H\big)(q\otimes 1)$$
is an $A\times_{\da,r}G/H-B\times_{\db,r}G/H$
imprimitivity bimodule. For this we only have to check that $A\times_{\da,r}G/H$ 
and $B\times_{\db,r}G/H$ are full corners in $L\times_{\dl,r}G/H$.
But since $p\otimes 1=\dl(p)$ it follows that
\begin{align*}
\big(L\times_{\dl,r}G/H\big)&(p\otimes 1)\big(L\times_{\dl,r}G/H\big)\\
&=\big((1\otimes M(C_0(G/H))\dl(L)\big)(p\otimes 1)
\big(\dl(L)(1\otimes M(C_0(G/H))\big)\\
&=(1\otimes M(C_0(G/H))\dl(LpL)(1\otimes M(C_0(G/H))
\end{align*}
which is dense in $L\times_{\dl,r}G/H$ because $LpL$ is dense in $L$.
The argument for $B\times_{\db,r}G/H$ is the same.
\end{proof}

\begin{proof}[Proof of Theorem \ref{thm-genimp1}]
Let $(X,\dx)$ be the Morita equivalence between 
${\delta}\,\widehat{\widehat{\ }}$ and $\delta$ of Proposition \ref{prop-Mor}.
Then Proposition \ref{prop-meq} provides a Morita equivalence
$X\times_{\dx,r}G/H$ between $B\times_{\delta,r}G/H$ and
$(B\times_{\delta}G\times_{\widehat{\delta},r}G)
\times_{{\delta}\,\widehat{\widehat{\ }},r}G/H$.
Now Green's imprimitivity theorem together with \cite{QS-RH}
provides a Morita equivalence $X_H^G$ between 
$(B\times_{\delta}G)\times_{\widehat{\delta},r}H$ and 
$(B\times_{\delta}G\times_{\widehat{\delta},r}G)
\times_{{\delta}\,\widehat{\widehat{\ }},r}G/H$.
Thus 
$$\widetilde{X}_H^G\otimes_{B\times G\times_r G\times_rG/H}(X\times_{\dx,r}G/H)$$
is a $(B\times_{\delta}G)\times_{\widehat{\delta},r}H-B\times_{\delta,r}G/H$
imprimitivity bimodule.
\end{proof}

\begin{rmk}
As we pointed out in the introduction, it would be preferable to have a more
concrete bimodule implementing the equivalence. We do not know whether the
original construction of Mansfield can be modified to avoid the assumption of
normality.
\end{rmk}

\end{document}